\documentclass[11pt]{amsart}
\usepackage{amscd,amssymb}
\theoremstyle{plain}
\newtheorem{Thm}[equation]{Theorem}
\newtheorem{Cor}[equation]{Corollary}
\newtheorem{Prop}[equation]{Proposition}
\newtheorem{Lem}[equation]{Lemma}

\numberwithin{equation}{section}

\DeclareMathOperator{\Image}{\mathrm{Im}}

\renewcommand{\O}{\mathcal{O}}
\newcommand{\la}{\langle\,}
\newcommand{\ra}{\,\rangle}
\renewcommand{\SS}{\mathcal{S}}
\newcommand{\C}{\mathbb{C}}
\newcommand{\g}{\mathfrak{g}}
\renewcommand{\a}{\mathfrak{m}_{\ell}}
\newcommand{\n}{\mathfrak{n}_{\ell}}
\newcommand{\m}{\mathfrak{m}}
\newcommand{\Ker}{\operatorname{Ker}}
\newcommand{\Ad}{\operatorname{Ad}}
\newcommand{\ad}{\operatorname{ad}}
\newcommand{\gr}{\operatorname{gr}}
\newcommand{\End}{\operatorname{End}}
\newcommand{\op}{\operatorname{op}}
\newcommand{\SL}{\operatorname{SL}}
\newcommand{\Wh}{\operatorname{Wh}}
\newcommand{\bplus}{{\mbox{$\bigoplus$}}}
\begin{document}
\title{QUANTIZATION OF SLODOWY SLICES}

\author{Wee Liang Gan}
\address{Department of Mathematics, University of Chicago, Chicago,
IL 60637, U.S.A.}
\email{wlgan@math.uchicago.edu}

\author{Victor Ginzburg}
\address{Department of Mathematics, University of Chicago, Chicago, 
IL 60637, U.S.A.}
\email{ginzburg@math.uchicago.edu}

\begin{abstract}
We give a direct proof of (a slight 
generalization of) the recent result of 
Premet related to generalized Gelfand-Graev
representations and of an equivalence  due to Skryabin.
\end{abstract}

\maketitle
\section{\bf Introduction}
\subsection{}
Let $\g$ be a complex semisimple Lie algebra, let $G$ be the adjoint group
of $\g$, let $\mathbb{O}$ be a nonzero nilpotent $\Ad G$-orbit in $\g$,
and let $e \in \mathbb{O}$. By the Jacobson-Morozov Theorem, there is an
$\mathfrak{sl}_{2}$-triple $(e,h,f)$ associated to $e$, i.e.
$[h,e]=2e$, $[h,f]=-2f$, $[e,f]=h$. Fix such an
$\mathfrak{sl}_{2}$-triple.
The \emph{Slodowy slice} to $\mathbb{O}$ at $e$ is defined to
be the affine space $e+\Ker\ad f$, see e.g. \cite{slodowy} \S 7.4.
The same slice has been used already by
 Harish-Chandra \cite{hc} \S12-14, see also
 \cite{ba}.

Since the Killing form on $\g$ is nondegenerate, there is an isomorphism 
${\Phi}:\g \rightarrow \g^{*}$ such that $\langle\Phi(e),f\rangle=1$.
Let $\chi = {\Phi}(e)$ and
$\SS = {\Phi}(e+\Ker\ad f)$. We will show in \S 
\ref{sec:poisson} that
$\SS$ has a natural Poisson structure. The aim of this paper is to
construct a quantization of $\SS$.

\subsection{}
Under the action of $\ad h$, we have a decomposition $\g = \bigoplus_{i
\in \mathbb{Z}} \g(i)$, where 
\[ 
\g(i) = \{ x \in \g \mid [h,x]=ix \}.
\] 
Note that there is a nondegenerate skew-symmetric bilinear form $\omega$ on
$\g(-1)$ defined by $\omega(x,y) = \chi([x,y])$ for all $x,y \in \g(-1)$. 
Fix an isotropic subspace $\ell$ of $\g(-1)$, and denote by 
$\ell^{{\perp}_{\omega}} 
\subset \g(-1)$ the annihilator of $\ell$ with respect to $\omega$.
Let 
\begin{displaymath}
\a = \ell\, \bplus\, \bigl(\oplus_{_{i \leq -2}}\; \g(i)\bigr) \qquad 
\textrm{and} \qquad
\n = \ell^{{\perp}_{\omega}}\, \bplus \,
\bigl(\oplus_{_{i \leq -2}}\; \g(i)\bigr).
\end{displaymath}
Note that $\a \subset \n$ and they are both nilpotent Lie subalgebras of
$\g$. Also, $\chi$ restricts to a character on $\a$.

Let $U\g$ and $U\a$ be the universal enveloping algebras of $\g$
and $\a$ respectively. Denote by $\C_{\chi}$ the 1-dimensional left
$U\a$-module obtained from the character $\chi$ of $\a$, and
let $Q_{\ell}=U\g \bigotimes_{U\a}
\C_{\chi}$ be the induced left $U\g$-module, equivalently,  the
quotient of $U\g$ by the left ideal $I_{\ell}$ generated by $x-\chi(x),$
for all $x \in \a$. Now consider the unique extension of the adjoint
action of $\n$ on $\g$ to derivations on $U\g$. Note that $I_{\ell}$ is
stable under this action of $\n$ because if $x \in \a$ and $n \in \n$,
then $(x-\chi(x))n = n(x-\chi(x)) + [x,n]$, and $\chi([x,n])=0$. Thus,
there is an induced $\ad \n$-action on $Q_{\ell}$. 

Let $H_{\ell} = Q_{\ell}^{\ad \n}$ be the subspace of all 
$x + I_{\ell} \in Q_{\ell}$ such that $nx-xn \in I_{\ell}$ for all
$n \in \n$.
Take any $x + I_{\ell}, y+I_{\ell} \in H_{\ell}$.
We define an algebra structure on $H_{\ell}$ by
$(x+I_{\ell})(y+I_{\ell})= xy+I_{\ell}.$
We claim that this multiplication is well defined. 
To see this, we note that for any $m\in\a$ and $y+I_{\ell} \in H_{\ell}$,
 by the definition
of  $H_{\ell}$ we have: $[m,y] \in [\a,y]\subset
[\n,y]\subset I_{\ell}$. Hence, we find: 
$(m-\chi(m))y= ym-y\chi(m)+[m,y] \in
y(m-\chi(m)) +I_{\ell}\subset I_{\ell}.$ 
 It follows that $I_{\ell}\cdot y\subset I_{\ell}$, and our claim is proved.
Further, 
it is  clear that   $H_{\ell}$ is
closed under the multiplication,
because if $n \in \n$, then
$nxy - xyn = (nx-xn)y + x(ny-yn) \in I_{\ell}$.

We will define in \S \ref{sec:kaz} the Kazhdan grading on $\C[\SS]$
and the Kazhdan filtration on $H_{\ell}$. We will then prove in \S 5
that 
$\gr H_{\ell}$ is canonically isomorphic to $\C[\SS]$ as graded Poisson
algebras, and $H_{\ell}$ is independent of the choice of
$\ell$. Our proof is based on
generalizations of some of Kostant's
results in \cite{kostant}, in which he considered the case when $e$ is a
principal nilpotent element.
Kostant's results were generalized by Lynch \cite{lynch} to the setup of
`admissible' parabolic subalgebras, which include our results in the special case of 
an even nilpotent  element $e$. In particular, our proof of Proposition \ref{thm:Q} below is very similar
to an argument in  \cite{lynch}, see also [Ko].
 However, we believe that an $\mathfrak{sl}_{2}$-triple
setting considered in the present paper
is more natural than that of admissible parabolic subalgebras considered
in \cite{lynch}, making all results much `cleaner'.

\subsection{}
This work was inspired by Premet \cite{premet}. He proved
the isomorphism between $\gr H_{\ell}$ and $\C[\SS]$ in the case when
${\ell}$ is Lagrangian.
His proof is based on results over algebraically closed fields of
positive characteristics. Similar results were also obtained in
\cite{tjin} using BRST cohomology.
Recall that in the Lagrangian case, $Q_{\ell}$ is called a
\emph{generalized Gelfand-Graev representation} associated to $e$ (c.f.
e.g.
\cite{kawanaka}, \cite{matumoto}, \cite{moeglin}, or \cite{yama}), we have
an algebra
isomorphism $\End_{U\g}(Q_{\ell})^{\op} \rightarrow H_{\ell}:h \to h(1
\otimes 1)$, and $H_{\ell}$ may be identified with the space of 
Whittaker vectors $\Wh (Q_{\ell})=\{ v \in Q_{\ell} \mid xv= \chi(x)v, \
\forall x \in \a \}$.

\section{\bf A decomposition lemma}
\subsection{}
It will be useful to define 
a linear action of $\C^{*}$ on $\g$ which stabilizes $e+\Ker \ad f$.
First, consider the Lie algebra homomorphism
$\mathfrak{sl}_{2}(\C) \rightarrow \g$ defined by
\begin{displaymath}
\left( \begin{array}{cc} 0 & 1 \\ 0 & 0 \end{array} \right) \mapsto e,
\quad
\left( \begin{array}{cc} 1 & 0 \\ 0 & -1 \end{array} \right) \mapsto h,
\quad
\left( \begin{array}{cc} 0 & 0 \\ 1 & 0 \end{array} \right) \mapsto f.
\end{displaymath}
This Lie algebra homomorphism exponentiates to a rational homomorphism
$\tilde{\gamma}: \SL_{2}(\C) \rightarrow G$. We put 
\begin{displaymath}
\gamma : \C^{*}
\rightarrow G\quad,\quad\gamma(t) = \tilde{\gamma}
\left( \begin{array}{cc} t & 0 \\ 0 & t^{-1} \end{array} \right)
\quad,\quad\forall t \in \C^{*}.
\end{displaymath}
Note that $\bigl(\Ad \gamma(t)\bigr)(e) = t^{2}\cdot e$.
The desired action of $\C^{*}$ on $\g$, to be denoted by
$\rho$, is defined by $\rho(t)(x)=t^{2}\cdot
\bigl(\Ad \gamma(t^{-1})\bigr)(x),$ for all
$x \in \g$. Note that $\rho(t)(e+x) = e + \rho(t)(x)$. Thus, since
$\rho(t)$ stabilizes $\Ker\ad f$, it also stabilizes $e+ \Ker\ad f$. Note
that 
$\underset{^{t\to 0}}{\lim}^{\,} \rho(t)(x) = e,
$ for all $x \in e+ \Ker\ad f$, i.e. the 
$\C^{*}$-action on $e+ \Ker\ad f$ is contracting.

\subsection{} \label{transversal}
Recall that the intersections of $e+ \Ker\ad f$ with $\Ad G$-orbits
in
$\g$ are transversal. This is clear at $e$ since $\g = [\g,e] \oplus (\Ker 
\ad f)$. Thus, the adjoint action map
$G \times (e+ \Ker\ad f) \rightarrow \g$ has a surjective differential at
each point
in some open neighborhood of $(1,e)$. It follows that at each point in
some open neighborhood of $e$ in $e+ \Ker\ad f$, the intersection of $e+
\Ker\ad f$ with $\Ad
G$-orbits is transversal. By the contracting $\C^{*}$-action on $e+
\Ker\ad f$, it
follows that the same is true at all points of $e+ \Ker\ad f$.

\subsection{}
Let $N_{\ell}$ be the unipotent subgroup of $G$ with Lie algebra $\n$, and
let $\a^{\perp} \subset \g^{*}$ be the annihilator of $\a$.
The following key lemma is a generalization of \cite{kostant} Theorem 1.2.

\begin{Lem} \label{isom} 
The coadjoint action map $\alpha : N_{\ell} \times \SS
\rightarrow \chi+ \a^{\perp}$ is an isomorphism of affine
varieties.
\end{Lem}
\begin{proof} Given a subspace $V\subset \g$,
we will write $V^{\perp_\g}$ for the annihilator
of $V$ in $\g$ (as opposed to
$V^\perp$, the annihilator in $\g^*$)
with respect to the Killing form on $\g$.
The statement of the Lemma 
is equivalent to saying that the adjoint action map
\begin{displaymath}
\alpha : N_{\ell} \times (e+ \Ker\ad f) \rightarrow e+
\a^{\perp_\g}
\end{displaymath}
is an isomorphism.
First, note that 
$\mathfrak{sl}_{2}$ representation theory implies:
 $[\n, e]\cap \Ker\ad f =0,$ and that the
map $\,\n\to [\n, e]\,,\, x\mapsto [x,e],$ is a bijection.
It follows easily that
$\dim\a^{\perp_\g}= \dim\n +\dim \g(0) +\dim\g(-1)
=\dim[\n, e]+\dim \Ker\ad f
$. Thus, there is a direct sum
decomposition
\begin{equation}\label{dec}
\a^{\perp_\g}=[\n, e]\,\oplus\, \Ker\ad f
\end{equation}

Next, define a $\C^{*}$-action on $N_{\ell} \times (e+ \Ker\ad f)$ by
\begin{displaymath}
t \cdot (g, x) = (\gamma(t^{-1})g\gamma(t), \rho(t)(x)).
\end{displaymath}
Note that for any $(g, x) \in N_{\ell} \times (e+ \Ker\ad f)$, we have
$\underset{^{t\to 0}}{\lim}^{\,}
t \cdot (g,x) = (1,e)$. Also, the action of $\C^{*}$ on $e +
\a^{\perp_\g}$ satisfies $\underset{^{t\to 0}}{\lim}^{\,}\rho(t)(x) = e,$ 
for any $x \in e + \a^{\perp_\g}$. The action map
$\alpha$ is  $\C^{*}$-equivariant. 
 Moreover, by (\ref{dec}), 
$\alpha$ induces an isomorphism between the tangent spaces of the
$\C^{*}$-fixed points $(1,e)$ and $e$. 
Thus, Lemma \ref{isom} follows from the following general result: 
{\it An equivariant
morphism $\alpha:X_{1} \rightarrow X_{2}$ of smooth
affine $\C^{*}$-varieties  with contracting $\C^{*}$-actions
which induces
an isomorphism between the tangent spaces of the
$\C^{*}$-fixed points must be an isomorphism.}

To prove this, let $x_{i}$ be the $\C^{*}$-fixed point of $X_{i}\,,\,
i=1,2$, let
$T_{i}$ be the tangent space of $X_{i}$ at $x_{i}$, and write
$\chi_{_{X}}
\in \C[[t]]$ for the formal character of the coordinate ring
of a $\C^{*}$-variety $X$. Since $\alpha$ induces
an isomorphism $T_{1} \stackrel{\sim}{\rightarrow} T_{2}$, the pullback
on coordinate rings $\alpha^{*}: \C[X_{2}] \rightarrow \C[X_{1}]$ is
injective. The surjectivity of $\alpha^{*}$ follows from
the equation: $\chi_{_{X_{1}}} = \chi_{_{T_{1}}} =\chi_{_{T_{2}}}
= \chi_{_{X_{2}}}$, see \cite{ginzburg}~(7.7).
\end{proof}

\medskip\noindent
{\bf Remark.}\, Even in the case of $e$ being the principal nilpotent,
the proof above is much simpler than that of \cite{kostant} Theorem 1.2.
There is also an alternative proof of Lemma \ref{isom} based on an inductive
argument similar to one used by Lynch \cite{lynch}.

\section{\bf Poisson structure on $\SS$} \label{sec:poisson}
\subsection{} \label{reduction0}
The space $\g^{*}$ has a natural Poisson structure defined by
\[ \{F_1,F_2\}(\xi) = \xi([dF_1(\xi), dF_2(\xi)]), \] 
where $F_1,F_2 \in \C[\g^{*}]$ and
$\xi \in \g^{*}$. The symplectic leaves of $\g^{*}$ are the
$\Ad^{*}G$-orbits. From \S \ref{transversal}, we know that $\SS$
intersects the symplectic leaves transversally. Thus,
to show that $\SS$ inherits a Poisson structure from
that of $\g^{*}$, it suffices (by \cite{vaisman} Proposition 3.10)
 to
verify that for any $\Ad^{*}G$-orbit $\mathcal{O}$ and $\xi \in
\mathcal{O} \cap \SS$, 
the restriction of the symplectic
form on $T_{\xi}\mathcal{O}$ to $T_{\xi}\SS \cap 
T_{\xi}\mathcal{O}= T_{\xi}(\SS \cap \mathcal{O})$ is nondegenerate. Note that 
$T_{\xi}\SS = {\Phi}(\Ker\ad f)$, and the annihilator of
${\Phi}(\Ker\ad f)$ in $\g$ is $[f, \g]$. 
Thus, the null space of
the restriction of the symplectic form to  $T_{\xi}\SS \cap
T_{\xi}\mathcal{O}$ is
\begin{displaymath}
{\Phi}\left(\big[{\Phi}^{-1}(\xi) \,,\,[f, \g]\big] \, \,\cap \, \, \Ker\ad f\right).
\end{displaymath}
This is $0$ since $\Phi^{-1}(\xi) \in (e+\Ker \ad f)$.

More details on the statements made above are provided in \S\ref{app}.

\subsection{} \label{reduction1}
The Poisson structure on $\SS$ can also be described via Hamiltonian
reduction. For this, we will take $\ell$ to be a Lagrangian subspace of
$\g(-1)$ in \S \ref{reduction1}. Then $\a =
\n$, and we denote both of them by $\m$. Let $M$ be
the unipotent subgroup of $G$ with Lie algebra $\m$. The moment
map $\mu : \g^{*} \rightarrow \m^{*}$ for the coadjoint action
of $M$ on $\g^{*}$  is just the restriction of functions from $\g$ to
$\m$.
Note that $\mu^{-1}(\chi|_{\m}) = \chi + \m^{\perp}$, where $\m^{\perp}
\subset \g^{*}$ is the annihilator of $\m$. 
Since $\chi|_{\m}$ is a character on $\m$, it is fixed under the
coadjoint action of $M$. 
Moreover,
$e+\m^{\perp_\g}$ is transversal to the $\Ad G$-orbits in
$\g$. This is clearly true at $e$, hence locally around $e$, and hence
it is true everywhere using the contracting $\C^{*}$-action $\rho$ (c.f.
\S \ref{transversal}). 
The transversality just proved implies in particular
that, for any $\xi\in \g$, we have $\g=[\g, \xi]+\m^{\perp_\g}$.
It follows that
$\chi|_{\m}$ is a regular value for
the restriction of $\mu$ to each  symplectic leaf of $\g^{*}$.
Thus, by Lemma \ref{isom}
 we have a Hamiltonian reduction of the Poisson structure on $\g^{*}$
to a
Poisson structure on $\SS$ (c.f. \cite{vaisman} Theorem 7.31).
We remark that the symplectic form on each symplectic leaf of $\SS$ is
obtained by symplectic reduction of the corresponding symplectic leaf of
$\g^{*}$. From this, and by the canonical embedding of $\SS$ into $\chi +
\m^{\perp}$, it is easy to see that the Poisson structures on $\SS$
defined in Sections \ref{reduction0} and \ref{reduction1} are the same. 
Moreover, the Poisson bracket $\{F_1, F_2\}$ for any $F_1, F_2 \in
\C[\SS]$ may be described as follows.
Let $\pi : (\chi +\m^{\perp}) \twoheadrightarrow
 (\chi +\m^{\perp})/\Ad^*M\simeq\SS$ be the projection map. Take
an arbitrary extension $\tilde{F_1}$ of $F_1
\circ \pi$ to $\g^{*}$, and an arbitrary extension $\tilde{F_2}$ of $F_2
\circ \pi$ to $\g^{*}$. Then $\{F_1, F_2\}\circ \pi =\{\tilde{F_1},
\tilde{F_2} \}\circ i$, where $i:
(\chi + \m^{\perp}) \hookrightarrow \g^{*}$.

\section{\bf Kazhdan grading and filtration} \label{sec:kaz}
\subsection{}
Now let us define a linear action $\rho^\sharp$ of $\C^{*}$ on $\g^{*}$
such that $\rho^\sharp$ stabilizes $\SS$. This action is defined by
$\rho^\sharp(t)(\xi)=t^{-2}\Ad^{*}(\gamma(t))(\xi)$ for all $\xi \in
\g^{*}$. 
If $\xi \in \g^{*}$ and $x\in \g$, then
$\langle\rho^\sharp(t)\xi, x\rangle
= t^{-4} \langle \xi, \rho(t)x\rangle$.

Note that we have an induced action on $\C[\g^{*}]$ defined by
$(\rho^\sharp(t)F)(\xi) = F(\rho^\sharp(t)^{-1}(\xi))$, where $F \in
\C[\g^{*}]$. 
The decomposition $\C[\g^{*}] = \bigoplus_{n \in \mathbb{Z}}
\C[\g^{*}](n)$,
where
\[
\C[\g^{*}](n) = \{F \in \C[\g^{*}] \mid \rho^\sharp(t)(F) = t^{n}F,
\quad
\forall t  \in \C^{*} \}, 
\]
gives $\C[\g^{*}]$ the structure of a graded
algebra. We call this the \emph{Kazhdan grading} on $\C[\g^{*}]$.

Similarly, since $\C^{*}$ acts on $\SS$ via $\rho^\sharp$, we may also
speak of the Kazhdan grading on $\C[\SS]$. Note however that the
weights of $\rho^\sharp$ on ${\Phi}(\Ker\ad f)$ are negative
integers. Thus, the
Kazhdan grading on $\C[\SS]$ has no negative graded components.

\subsection{}
Let $S\g$ be the symmetric algebra of $\g$. Identify $S\g$ with
$\C[\g^{*}]$. From this identification, $S\g$ acquires a graded algebra
structure from the Kazhdan grading on $\C[\g^{*}]$. This grading on $S\g$
may
be described explicitly as follows. Let $S\g = \bigoplus_{n \geq 0}
S^{n}\g$ be the standard grading of $S\g$. The action of $\ad h$ on $\g$
extends uniquely to a derivation on $S\g$. For any $i \in \mathbb{Z}$,
let 
\[
(S^{n}\g)(i) =\{ x \in S^{n}\g \mid (\ad h)(x) = ix \}.
\]
The Kazhdan grading $S\g = \bigoplus_{n \in \mathbb{Z}} (S\g)[n]$ is
defined by letting $(S\g)[n]$ be the subspace of $S\g$ spanned by all 
$(S^{j}\g)(i)$ with $i+2j=n$.

We will define a
filtration on $U\g$ such that the associated graded algebra is $S\g$ with
the Kazhdan grading.
Consider the standard filtration $U_{0}\g \subset U_{1}\g \subset \ldots
\subset U_{n}\g \subset \ldots$ \ of $U\g$. The action of of $\ad h$ on
$\g$ extends
uniquely to a derivation on $U\g$ which we will also denote by $\ad h$.
For any $i \in \mathbb{Z}$, let 
\[
(U_{n}\g)(i) = \{x \in U_{n}\g \mid (\ad h)(x) = ix\}.
\]
The \emph{Kazhdan
filtration}, $\ldots \subset F_{n}U\g
\subset F_{n+1}U\g \subset \ldots,$ \ of $U\g$ is a
$\mathbb{Z}$-filtration defined by letting
$F_{n}U\g$ be the subspace of $U\g$ spanned by all $(U_{j}\g)(i)$ with
$i+2j \leq n$. The Kazhdan filtration gives $U\g$ the
structure of a filtered algebra. Observe that if $x \in F_{n}U\g$ and $y \in
F_{m}U\g$, then $xy-yx \in F_{n+m-2}U\g$. Further,
for any $x\in \g(n)\,,\,y\in \g(m)$, we have
$x\in F_{n+2}U\g,$\
$y\in F_{m+2}U\g$ and, moreover, 
the classes of $x$ and $y$ in $\gr U\g$ commute.
It follows that the natural map $\g \to \gr U\g$ extends uniquely
to a well-defined graded algebra homomorphism: $S\g\to \gr U\g$.
Moreover, the Poincar\'{e}-Birkhoff-Witt Theorem implies easily that this homomorphism
is a bijection.

Note further that the canonical Poisson bracket on $\gr U\g$
 defined by
$\{\gr_{n}x, \gr_{m}y \}= \gr_{n+m-2}(xy-yx)$
corresponds, under the
 isomorphisms $\gr U\g= S\g=\C[\g^{*}]$, to
the standard  Poisson structure on
$\C[\g^{*}]$ (it suffices to check this on linear functions).

\subsection{}
Let $p:U\g \rightarrow Q_{\ell}$ be the quotient map.
The Kazhdan filtration on $Q_{\ell}$ is defined by
$F_{n}Q_{\ell}=p(F_{n}U\g)$.
This
gives $Q_{\ell}$ the structure of a filtered module over $U\g$. Note that
$\gr Q_{\ell}$
has a commutative 
 algebra structure since $\gr Q_{\ell}=\gr U\g / \gr I_{\ell}$,
where $I_{\ell}$ is the left ideal generated by $x -\chi(x),$ for all $x
\in \a$. Further, $F_{n}Q_{\ell} = 0$ for all $n<0$.
The associated graded map $\gr p: S\g
\rightarrow \gr Q_{\ell}$ is a surjective homomorphism of graded algebras
whose
kernel is the ideal generated by $x-\chi(x)$ for all $x \in \a$. Under the
identification of $S\g$ with $\C[\g^{*}]$, we see that $\Ker \gr p$
is the ideal of all polynomial functions on $\g^{*}$ vanishing on 
$\chi + \a^{\perp}$.
Since $\a^{\perp}$ is stable under the action of $\C^{*}$
via $\rho^\sharp$, there is a Kazhdan grading on
$\C[\chi +\a^{\perp}]$, and
$\gr Q_{\ell}$
may be identified with $\C[\chi+\a^{\perp}]$ as graded algebras. Note
that
the weights of $\rho^\sharp$ on $\a^{\perp}$ are negative integers, which
agrees
with the fact that the Kazhdan grading on $\gr Q_{\ell}$ has no negative
component.

\subsection{}
The Kazhdan filtration on $H_{\ell}$ is induced from the Kazhdan
filtration on $Q_{\ell}$ via the inclusion $H_{\ell} \hookrightarrow
Q_{\ell}$.
Note that $\gr H_{\ell} \hookrightarrow \gr Q_{\ell}$ is an injective
homomorphism of 
graded algebras. 
Since ${\Phi}(\Ker \ad f) \subset \a^{\perp}$, we have a restriction
homomorphism $\nu: \C[\chi+\a^{\perp}] \rightarrow \C[\SS]$. Recall 
that $\gr Q_{\ell}$ is identified with $\C[\chi+\a^{\perp}]$. Thus,
there
is a canonical homomorphism of graded algebras 
$\nu: \gr H_{\ell} \rightarrow \C[\SS]$.

Our goal is to give a simple direct proof of the following theorem.
\begin{Thm} \label{thm:premet} 
The canonical homomorphism $\nu : \gr H_{\ell} \rightarrow \C[\SS]$
is an isomorphism of graded Poisson algebras. Moreover, $H_{\ell}$ is
independent of the choice of an isotropic subspace $\ell
\subset~\g(-1).$
\end{Thm}

\section{\bf $\n$-cohomology of $\gr Q_{\ell}$ and $Q_{\ell}$}
\subsection{}
From now on, we will
regard  $U\g$ and 
$Q_{\ell}$ as a $\n$-modules via the adjoint $\n$-action.
Thus, $H_{\ell} = H^{0}(\n\,,\, Q_{\ell})$. 
Note that the map $p:U\g \rightarrow
Q_{\ell}$ is $\n$-equivariant. 

Next, note that $\n$ is a graded subalgebra of $\g$, and hence it is also
filtered. Clearly, $U\g$ and $Q_{\ell}$ are Kazhdan filtered $\n$-modules.
Thus,
$\gr U\g$ and $\gr Q_{\ell}$ acquire the
structure of Kazhdan graded $\n$-modules, and $\gr p: \gr U\g \rightarrow 
\gr Q_{\ell}$ is $\n$-equivariant. The claim that $\nu : \gr H_{\ell}
\rightarrow \C[\SS]$ is an isomorphism 
follows immediately from the following two propositions.

\begin{Prop} \label{thm:grQ}
$\nu: H^{0}(\n\,,\, \gr Q_{\ell}) \rightarrow \C[\SS]$ is an
isomorphism,
and $H^{i}(\n\,,\, \gr Q_{\ell}) = 0$ for all $i>0$.
\end{Prop}

\begin{Prop} \label{thm:Q}
$\gr H^{0}(\n\,,\, Q_{\ell}) = H^{0}(\n\,,\, \gr Q_{\ell})$, and $H^{i}(\n\,,\,
Q_{\ell}) = 0$ for all $i>0$.
\end{Prop}

\subsection{}
To prove the above propositions, let us describe an action of $\n$ on
$\C[\chi+ \a^{\perp}]$ so that the identification of $\gr Q_{\ell}$
with $\C[\chi+ \a^{\perp}]$ becomes also an identification of Kazhdan
graded $\n$-modules.

First, note that the
adjoint action of $\n$ on $\g$ extends uniquely to derivations on $S\g$.
We have $\gr U\g = S\g$ as Kazhdan graded $\ad \n$-modules.
Note also that the $\Ad^{*} N_{\ell}$-action on $\g^{*}$ induces an
$N_{\ell}$-action on $\C[\g^{*}]$.
This $N_{\ell}$-action on $\C[\g^{*}]$ is
locally finite, so there is an infinitesimally induced action of
$\n$ on $\C[\g^{*}]$. We have $S\g=\C[\g^{*}]$ as Kazhdan graded $\ad
\n$-modules.
Since $\chi+\a^{\perp}$ is stable under the
$\Ad^{*}N_{\ell}$-action, we also have an
induced $N_{\ell}$-action on $\C[\chi+\a^{\perp}]$, and hence also
an
infinitesimal  action of $\n$ on $\C[\chi+\a^{\perp}]$.
Clearly,
the restriction homomorphism $\C[\g^{*}] \rightarrow
\C[\chi+\a^{\perp}]$ is
$\n$-equivariant. Hence, $\gr Q_{\ell}$ and $\C[\chi+\a^{\perp}]$ are
identified as $\n$-modules.

\subsection{}
We now prove Proposition \ref{thm:grQ}. 
Note that
if we let $N_{\ell}$ act on $N_{\ell} \times \SS$ by left
translations on
$N_{\ell}$, then $\alpha$ is $N_{\ell}$-equivariant. We have an
$N_{\ell}$-module
structure on $\C[N_{\ell}]$ defined by $(g\cdot F)(x) = F(g^{-1}x)$,
where
$g,x \in N_{\ell}$ and $F \in \C[N_{\ell}]$.
Lemma \ref{isom} implies
that $\C[N_{\ell}] \otimes \C[\SS] \cong \C[\chi+ \a^{\perp}]$ as
algebras and
$N_{\ell}$-modules. Since $H^{0}(\n\,,\, \C[\chi+ \a^{\perp}])$ is
precisely the
$N_{\ell}$-invariants $\C[\chi+ \a^{\perp}]^{N_{\ell}}$, it follows
that
\begin{displaymath}
H^{0}(\n\,,\, \gr Q_{\ell}) = H^{0}(\n\,,\, \C[\chi+ \a^{\perp}]) = 
\C[\chi+ \a^{\perp}]^{N_{\ell}}
\stackrel{\nu}{\stackrel{\sim}{\longrightarrow}} \C[\SS].
\end{displaymath}
Also, for $i > 0$,
\begin{displaymath}
H^{i}(\n\,,\, \gr Q_{\ell}) = H^{i}(\n\,,\, \C[\chi+ \a^{\perp}])
= H^{i}(\n\,,\, \C[N_{\ell}]) \otimes \C[\SS] = 0,
\end{displaymath}
where the last equality follows from the fact that the standard cochain
complex for Lie algebra cohomology with coefficients in $\C[N_{\ell}]$ is
just the algebraic
de Rham complex for $N_{\ell}$, and $N_{\ell}$ is isomorphic to
an affine space, hence has trivial de Rham cohomology.
This completes the proof of Proposition \ref{thm:grQ}.

\subsection{}
We now deduce Proposition \ref{thm:Q} from Proposition \ref{thm:grQ} using 
spectral sequence. Note that $\n$ is a negatively graded subalgebra of
$\g$, so its dual $\n^{*}$ is positively graded; we write its
decomposition as $\n^{*} = \bigoplus_{i \geq 1}\n^{*}(i)$.
Consider the standard
cochain complex for computing the $\n$-cohomology of $Q_{\ell}$:
\begin{equation} \label{eqn:Q}
0 \longrightarrow Q_{\ell} \longrightarrow \n^{*} \otimes Q_{\ell}
\longrightarrow
\ldots \longrightarrow \wedge^{n}\n^{*} \otimes Q_{\ell} \longrightarrow
\ldots.
\end{equation}
A filtration on $\wedge^{n}\n^{*} \otimes Q_{\ell}$ is defined by letting
$F_{p}(\wedge^{n}\n^{*} \otimes Q_{\ell})$ be the subspace of
$\wedge^{n}\n^{*} \otimes
Q_{\ell}$ spanned by $(x_{1} \wedge \ldots \wedge x_{n}) \otimes v$, for
all
$x_{1} \in \n^{*}(i_{1}), \ldots, x_{n} \in \n^{*}(i_{n})$ and $v \in
F_{j}Q_{\ell}$
such that $i_{1} + \ldots + i_{n} + j \leq p$. This defines the structure
of a filtered complex on (\ref{eqn:Q}). Taking the associated graded gives
the standard cochain complex for computing the $\n$-cohomology of $\gr
Q_{\ell}$.

Now consider the spectral sequence with
\begin{displaymath}
E_{0}^{p,q} = F_{p}(\wedge^{p+q}\n^{*} \otimes Q_{\ell})/
F_{p-1}(\wedge^{p+q}\n^{*}
\otimes Q_{\ell}).
\end{displaymath}
Then $E_{1}^{p,q} = H^{p+q}(\n\,,\, \gr_{p} Q_{\ell})$.
Hence, Proposition \ref{thm:Q} follows from Proposition \ref{thm:grQ} and
the fact
that the spectral sequence converges to
$E_{\infty}^{p,q} = F_{p}H^{p+q}(\n\,,\, Q_{\ell})/ F_{p-1}H^{p+q}(\n\,,\,
Q_{\ell})$.

\subsection{}
Note that if $\ell_{1} \subset \ell_{2}$ are both isotropic subspaces of
$\g(-1)$, then we have a natural map $Q_{\ell_{1}} \rightarrow
Q_{\ell_{2}}$ which gives a map $H_{\ell_{1}} \rightarrow H_{\ell_{2}}$.
By above, the associated graded map $\gr H_{\ell_{1}}
\rightarrow \gr H_{\ell_{2}}$ is an isomorphism of Kazhdan graded
algebras,
hence $H_{\ell_{1}} \rightarrow H_{\ell_{2}}$ is an isomorphism of Kazhdan
filtered algebras.
Taking $\ell_{1} = 0$, we see that $H_{\ell_{2}}$ is independent of the
choice of $\ell_{2}$.

To check that the Poisson structure of $H_{\ell_{1}}$ is that of
$\C[\SS]$, we take $\ell_{2}$ to be a Lagrangian subspace. It suffices
to prove that $\gr H_{\ell_{2}}$ and $\C[\SS]$ have the same Poisson
structures. This follows from the observation that the Poisson structure
of the Kazhdan graded algebra $\gr U\g$ is same as the usual Poisson
structure of $\C[\g^{*}]$ and the discussion in \S \ref{reduction1}.

This concludes the proof of Theorem \ref{thm:premet}.

\section{\bf Skryabin's equivalence}
\subsection{}
We would like to make some remarks on  results by Skryabin [Sk].
Throughout this section we
take $\ell$ to be a Lagrangian subspace of $\g(-1)$, and write
$\m = \a$, $Q = Q_{\ell}$, and $H = H_{\ell}$. Let $\mathcal{C}$ be the
abelian
category of finitely generated
$U\g$-modules on which $m-\chi(m)$ acts locally nilpotently for
each $m \in \m$. For $E \in \mathcal{C}$, put $\Wh(E) =
\{ x \in E \mid$
$ mx=\chi(m)x\,,\,\forall m \in \m \}$. 
Observe that for $E \in \mathcal{C}$ we have:
$\Wh(E) =0\enspace\Longrightarrow\enspace E=0.$ 

The following beautiful theorem of Skryabin and its proof
in [Sk] are similar in spirit to the well-known 
result of Kashiwara on the equivalence of the
category of ${\mathcal{D}}$-modules on a submanifold with
the category of ${\mathcal{D}}$-modules on the ambient manifold
which are supported on the submanifold.

\begin{Thm} \label{skr}
The functor $V \mapsto Q \otimes_{H} V$ sets up an equivalence
of the category of finitely generated left $H$-modules
and category $\mathcal{C}$. The inverse  equivalence
is given by the functor $E \mapsto \Wh(E)$, in particular,
the latter is exact.
\end{Thm}

\subsection{} We present an alternative proof of Theorem \ref{skr}
along the lines of the preceeding section
(cf. also [Ko] \S4, and [Ly] Theorem 4.1).

Fix an $H$-module $V$ generated by a finite dimensional subspace $V_0$.
View $H$ as a filtered algebra with respect to the Kazhdan filtration
$F_\bullet H$, and define an increasing filtration on $V$
by $F_iV= F_iH\cdot V_0$, for all $i$. This makes
$V$ a filtered $H$-module, with associated graded
$\gr H$-module $\gr V$. Lemma \ref{isom} and Theorem 
\ref{thm:premet} yield
$\gr Q =
\C[N]\otimes \gr H$. We deduce
\begin{displaymath}
H^{0}\bigl(\m\,,\, \gr Q \otimes_{_{\gr H}} \gr V\bigr)
=\gr V\quad\mbox{and}\quad
H^{i}\bigl(\m\,,\, \gr Q \otimes_{_{\gr H}} \gr V\bigr)=0\;,\;
\forall i>0.
\end{displaymath}

Further, the  Kazhdan filtration on  $Q$
and the filtration on $V$ give rise to a filtration on
$Q \otimes_{H} V$, and since 
$\gr Q$ is free over $\gr H$,
we have
 a canonical isomorphism
$\,\gr\bigl(Q \otimes_{H} V\bigr)\,\simeq\,
\gr Q \otimes_{\gr H} \gr V\,.$
Now a spectral sequence argument very similar to that in the proof
of Proposition \ref{thm:Q} yields 
\begin{equation}\label{Wh-vanish}
H^{0}\bigl(\m\,,\, Q \otimes_{H} V\bigr)
= V\quad\mbox{and}\quad
H^{i}\bigl(\m\,,\,Q \otimes_{H} V\bigr)=0\;,\;
\forall i>0\,,
\end{equation}
where the cohomology is taken with respect to the
$\chi$-twisted action of $\m$.

Note that the equation on the left of (\ref{Wh-vanish}) says that
$\Wh(Q \otimes_{H} V)=V$.
Thus, to complete the proof of the Theorem it suffices to show that,
for any $E\in {\mathcal{C}}$, the canonical map
$f: Q \otimes_{H} \Wh(E) \to E$ is an isomorphism.
Let $E'$ be the kernel, and $E''$ the cokernel of $f$.
We observe that $\Wh(E')=
E'\cap \Wh\bigl(Q \otimes_{H} \Wh(E)\bigr)$, which is
equal to $E'\cap \Wh(E)$, by (\ref{Wh-vanish}).
But $\Wh(E)$ does not intersect the kernel of $f$, hence
$\Wh(E')=E'\cap \Wh(E)=0$. Since $E'$ is clearly an object of
${\mathcal{C}}$, this yields $E'=0$. Hence $f$ is injective.

To prove surjectivity we write the long exact sequence of
cohomology associated to the short exact sequence
$\,0\to Q \otimes_{H} \Wh(E) \to E\to E''\to 0$. We obtain:
$$
0\to
H^{0}\bigl(\m\,,\, Q \otimes_{H} \Wh(E)\bigr)
\,\stackrel{H^0(f)}{\longrightarrow}\,
H^{0}(\m, E)\to H^{0}(\m, E'') \to
H^1\bigl(\m\,,\, Q \otimes_{H} \Wh(E)\bigr)\to\ldots\,.
$$
In this formula, 
the $H^1$-term vanishes,
and the map $H^0(f)$ is a bijection, due to (\ref{Wh-vanish}).
Hence, the  long exact sequence yields
$\Wh(E'')=H^{0}(\m, E'')=0$.
This forces $E''=0$, and the Theorem is proved.

\subsection{} 
Let $M$ be the unipotent algebraic subgroup of $G$
corresponding to the Lie algebra $\m$, and let $\mathcal{B}$ denote
the flag manifold for $G$.
Let  $\mathcal{V}$ be  an
 $M$-equivariant
$\mathcal{D}$-module on $\mathcal{B}$.
Given an element  $x\in \m$, we write 
$x_{_\mathcal{D}}$ for the action on $\mathcal{V}$
of the vector field corresponding
to $x$ via the $\mathcal{D}$-module structure,
and $x_{_M}$ for the action on $\mathcal{V}$
obtained by differentiating the $M$-action
arising from the equivariant structure.
We say that  $\mathcal{V}$ 
is an $\m$-Whittaker $\mathcal{D}$-module with respect to the character
$\chi: \m \rightarrow \C$ if, for any $x\in \m$ and $v\in \mathcal{V}$, we have
$(x_{_\mathcal{D}}-x_{_M})v=\chi(x)\cdot v$.

Note that the natural map $U\g \rightarrow Q$ maps the center of $U\g$
injectively into $H$ (c.f. \cite{premet} \S 6.2).
Denote by $Z_{+}$ the augmentation ideal of the center of $U\g$. 
Skryabin's result combined with the Beilinson-Bernstein localization
theorem implies the following.
\begin{Prop}
The category of finitely generated
$H/Z_{+}H$-modules is equivalent to the category of
$\m$-Whittaker coherent $\mathcal{D}_{\mathcal{B}}$-modules 
(with respect to the character
$\chi$).
\end{Prop}

\section{\bf Appendix:\ Intersections of a Slodowy slice with coadjoint orbits}
\label{app}
In this appendix, we provide more details on the
arguments in \S 3.1. These details were omitted in
the original version of the paper published more than ten years ago;
we had thought that those were routine arguments but since then we have
not always been able to recall those arguments ourselves.

Recall the setup and notation of \S\ref{reduction1}.
Restricting  the map $\mu$ to a coadjoint orbit   ${\mathcal O}$
gives a moment map associated with the $M$-action
on
${\mathcal O}$. According to  \S\ref{reduction1},
the element  $\chi|_{\mathfrak m}\in {\mathfrak m}^*$ is a regular
value of the map $\mu|_{\mathcal O}$. 
Thus,
by standard results concerning  fibers of  moment maps 
over regular values, cf. eg. \cite[Theorem 2.5]{GS},
we obtain the following

\begin{Lem}\label{nil} Assume that 
the set $\Sigma:=\O\cap (\chi +{\mathfrak m}^\perp)$ is non-empty.

Then, $\Sigma$ is a smooth coisotropic submanifold of the symplectic manifold
$\O$. Moreover, the canonical
nil-foliation on the coisotropic manifold $\Sigma$ is
the foliation by the $M$-orbits in  $\Sigma$.
\end{Lem}

Now, fix $\xi\in{\mathcal O}\cap{\mathcal S}$ and
let $T_\xi (M\cdot \xi)$ be the tangent space
to the $M$-orbit of $\xi$. By Lemma \ref{nil}, 
$T_\xi( M\cdot \xi)$ is the tangent space to the
nil-leaf through $\xi$. Further, the isomorphism
of Lemma 2.1 implies a direct sum decomposition
\[{\mathfrak m}^\perp=T_\xi (M\cdot\xi)\ \oplus\  T_\xi{\mathcal S}=
\ad{\mathfrak m} (\xi)\ \oplus\  \Phi(\ker\ad f).\]

Now, the claim from  \S 3.1 is an immediate consequence
of the above results, that is, 
we get

\begin{Cor}\label{non}
The manifold ${\mathcal O}\cap{\mathcal S}$ is 
a  symplectic submanifold of  ${\mathcal O}$,
for any coadjoint orbit ${\mathcal O}$;
that is, the
Kirillov-Kostant $2$-form on
${\mathcal O}$ gives, by restriction,
a nondegenerate  $2$-form on  ${\mathcal O}\cap{\mathcal S}$.
\end{Cor}

Here is an alternative, more direct proof of 
Corollary \ref{non}.

We use the notation of \S3.1 and let $x=\Phi^{-1}(\xi)$.
We identify $[x,\g]$ with $T_\xi\O$,
 write $\la-,-\ra$ for the Killing form on $\g$
and $(-)^{\perp_\g}\subset \g$ for the annihilator
with respect to  the Killing form.
The Kirillov-Kostant symplectic form on $[x,\g]$ is then given by the formula
$$[x,\g] \times [x,\g] \longrightarrow \mathbb C, \quad
[x,u]\times [x,v] \mapsto \la x,[u,v] \ra.$$

Now let $y=[x,v]\in T_\xi\O$ be such  that
$\la x,[u,v] \ra = 0$ holds
for all $u\in \g$ such that $[x,u] \in \Ker(\ad(f))$. 
Our formula (3.2) claims that
\begin{equation}\label{32}
y \in [x, [f,\g]]. 
\end{equation}

To see this, note first that  $\la x,[u,v]\ra = \la [x,u],v\ra$, so
$$\la a, v \ra = 0 \quad \mbox{ for all } a\in \Image(\ad x)\cap \Ker(\ad f).$$
Thus, our assumption on $y$ reads
$$v\in \Big([x,\g]\cap \Ker(\ad f) \Big)^{\perp_\g}. $$
For any $a\in\g$, the linear map $\ad a:\g\to\g$  
is skew-adjoint relative to the Killing form.
Therefore, we obtain
$$
  \big(\Image(\ad x)\cap \Ker(\ad f) \big)^{\perp_\g} = \Image(\ad x)^{\perp_\g} + \Ker(\ad f)^{\perp_\g} 
 = \Ker(\ad x) + \Image(\ad f).
$$

This proves \eqref{32} since we have
\[ y=[x,v] \in [x, \Ker(\ad x)]+[x, [f,\g]]=[x, [f,\g]].\]

\medskip
To complete the proof we must show that, for any  $x\in
e+\Ker(\ad f)$,
one has
\begin{equation}\label{ef}
[x, [f,\g]]\ \cap\  \Ker(\ad f)=0.
\end{equation}
To this end, we
use  direct sum decompositions
\begin{equation}\label{dir}
[e,\g]\ \oplus\  \Ker(\ad f)\ =\ \g\ =\ [f,\g]\ \oplus\ \Ker(\ad e).
\end{equation}

From the 
decomposition on the right we see that the map
$\ad e:\ [f,\g]\to [e,[f,\g]]$ is a bijection.
We deduce by continuity that, 
for any
$x\in \g$ sufficiently close to $e$, 
the map
$\ad x:\ [f,\g]\to [x,[f,\g]]$ is also a bijection.
It follows, that assigning to $x$ the
vector space $[x,[f,\g]]$ gives a continuous map
of a neighborhood of $e$ in $\g$ to an appropriate
Grassmannian.

Next, we observe that the first decomposition  in \eqref{dir} implies that
 $[e, [f,\g]] \cap \Ker(\ad f)=0$.
Therefore,  \eqref{ef} holds
for all $x$ in a neighborhood of $e$,
by continuity. This implies \eqref{ef}
for any  $x\in e+\Ker(\ad f)$,  using the ${\mathbb C}^*$-action.

There is also a completely algebraic proof of  \eqref{ef} as follows.
Let $y\in [x, [f,\g]] \cap \Ker(\ad(f))$, where $x=e+k$ for some $k\in
\Ker(\ad f)$.
Then, for some $z\in \g$, 
$$ y= [e+k,[f,z]] = [e,[f,z]] + [k,[f,z]] .$$
Suppose that $[f,z]\neq 0$. If $z\in \g(i)$,
then $[e,[f,z]]$ is a nonzero vector contained in $\g(i)$,
but $[k,[f,z]]$ is contained in $\oplus_{j<i} \g(j)$. 
More generally, we can write $$z=z_{i_1}+\cdots+z_{i_r}, \quad z_i\in \g(i) ,
\quad  i_1<\cdots<i_r.$$
Suppose $i_s$ is the biggest $i$ such that $[f,z_i]\neq 0$.
Then $[e,[f,z_{i_s}]]$ is a nonzero vector in $\g(i_s)$, but
$[e,[f,z_{i_j}]]$ for $j\neq s$ and $[k,[f,z]]$ are all contained in $\oplus_{j<i_s} \g(j)$. 
Therefore, if we decompose $y$ 
according to the direct sum $\g = \Ker(\ad f) \oplus [e,\g]$, then $y$ will have a nonzero
component in $[e,\g]$, which is a contradiction to $y\in \Ker(\ad f)$. Hence, 
$[f,z]=0$.

\bibliographystyle{plain}

\end{document}